\documentclass[11pt,a4paper]{article}
\usepackage[english]{babel}

\usepackage[dvips]{graphicx}
\usepackage[svgnames]{xcolor}

\colorlet{foldline}{Red}
\colorlet{crease}{Gray}
\colorlet{halfgray}{DarkGrey}
\colorlet{frontside}{GhostWhite!100!WhiteSmoke}
\colorlet{backside}{LightSteelBlue!25!LightGray}

\usepackage{amssymb}
\usepackage[tbtags]{amsmath}
\usepackage{amsthm}
\usepackage{icomma}

\renewcommand{\Re}{\operatorname{Re}}

\usepackage[T1]{fontenc}
\usepackage[osf,sc]{mathpazo}
\usepackage[euler-digits,euler-hat-accent]{eulervm}

\usepackage[top=32mm]{geometry}

\usepackage{caption}
\usepackage{fancyhdr}
\usepackage{hyperref}
\hypersetup{
  breaklinks = true,
  linktocpage=true,
  linkcolor  = MidnightBlue,
  citecolor  = MidnightBlue,
  urlcolor   = MidnightBlue,
  colorlinks = true,
}
\usepackage{breakurl}
\usepackage{breakcites}    

\title{Construction of a regular hendecagon by two-fold origami}
\author{Jorge C. Lucero\thanks{Dept.\ Computer Science, University of Bras\'{i}lia, Brazil. E-mail: \href{mailto:lucero@unb.br}{lucero@unb.br} }} 
\date{\today} 

\begin{document}

\maketitle

\begin{abstract} 
\noindent 
The regular hendecagon is the polygon with the smallest number of sides that cannot be constructed by single-fold operations of origami on a square sheet of paper. This article shows its construction by using an operation that requires two simultaneous folds.
\end{abstract}

\section{Introduction}

Single-fold origami refers to geometric constructions on a sheet of paper by performing a sequence of single folds, one at a time \cite{Alperin2006}. Each folding operation achieves a minimal set of specific incidences (alignments) between given points and lines by folding along a straight line, and there is a total of eight possible operations \cite{Lucero2017}.
The set of single-fold operations allows for the geometric solution of arbitrary cubic equations \cite{Alperin2000,Geretschlager1995}. As a consequence, the operations may be applied to construct regular polygons with a number of sides $n$ of the form $n=2^r3^sp_1p_2\ldots p_k$, where $r, s, k$ are nonnegative integers and $p_1, p_2,\ldots, p_k$ are distinct Pierpont primes of the form $2^m3^n+1$, where $m, n$ are nonnegative integers \cite{Gleason1988}. For example, previous articles in \textit{Crux Mathematicorum} have shown the construction of the regular heptagon \cite{Geretschlaeger1997} and nonagon \cite{Geretschlager1997b}. Let us note that this family of regular polygons is the same that can be constructed by straightedge, compass and angle trisector. 

Number 11 is the smallest integer not of the above form (the next are 22, 23, 25, 29, 31,\ldots); therefore, the hendecagon is the polygon with the smallest number of sides that can not be constructed by single-fold origami. In fact, its construction requires the solution of a quintic equation \cite{Gleason1988}, which can not be obtained by single folds. It has been shown that any polynomial equation of degree $n$ with real solutions may be geometrically solved by performing $n-2$ simultaneous folds \cite{Alperin2006}. Hence, a quintic equation could be solved by performing three simultaneous folds. However, a recent paper presented an algorithm for solving arbitrary quintic equations with only two simultaneous folds \cite{Nishimura2015}.  In this article, the algorithm will be applied to solve the quintic equation associated to the regular hendecagon, and full folding instructions for its geometric construction will be given.

Let us note that another problem related to a quintic equation, the quintisection of an arbitrary angle, has also been solved by two-fold origami \cite{Lang2004}.

\section{\label{S1}Quintic equation for the regular hendecagon and its solution}

\subsection{Equation}

A similar approach to that used for the heptagon \cite{Geretschlaeger1997} is followed. 

\begin{figure}
\centering
\includegraphics{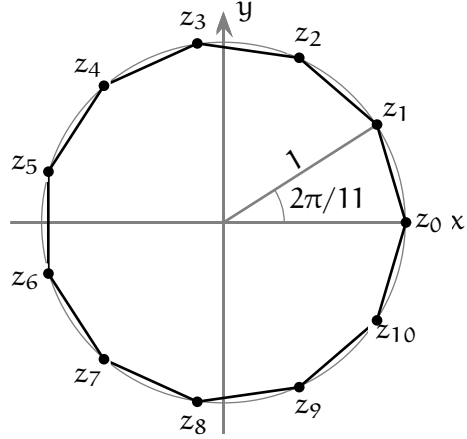}
\caption{The regular hendecagon.} 
\label{hend}
\end{figure}

Consider an hendecagon in the complex plane, inscribed in a circle of unitary radius (Fig.~\ref{hend}). Its vertices are solutions of the equation
\begin{equation}
z^{11}-1=0
\end{equation} 
One solution is $z=z_0=1$, and the others are solutions of
\begin{equation}
\frac{z^{11}-1}{z-1}=z^{10}+z^{9}+z^{8}+\cdots +z+1=0
\end{equation} 
Assuming $z\neq 0$ and dividing both sides by $z^5$ produces
\begin{equation}
z^{5}+z^{4}+z^{3}+\cdots +\frac{1}{z^3}+\frac{1}{z^4}+\frac{1}{z^5}=0
\label{eq3}
\end{equation}

Next, let $t=z+\overline{z}=2\Re z$. Since $|z|=|\overline{z}|=1$, then $\overline{z}=z^{-1}$ and $t=z+z^{-1}$. Expressing Eq.~(\ref{eq3}) in terms of $t$, the equation reduces to
\begin{equation}
t^{5}+t^{4}-4t^{3}-3t^2+3t+1=0
\label{eqt}
\end{equation} 

Solutions of Eq.~(\ref{eqt}) have the form $t_k= 2\Re z_k=2\cos(2k\pi/11)$, where $k=1, 2,\ldots, 5$ and $z_k$ are vertices as indicated in Fig.~\ref{hend}. Note that, due to symmetry on the $x$-axis, vertices except $z_0$ appear as complex conjugate pairs with a common real part. 

\subsection{Solution}

\begin{figure}[!htb]
\centering
\includegraphics{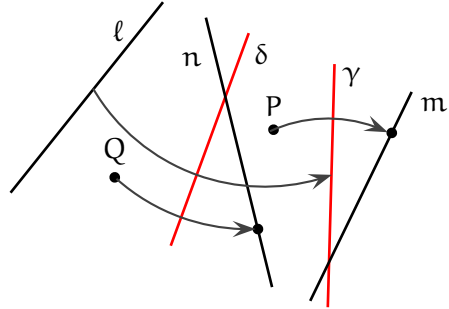}
\caption{A two-fold operation. Red lines $\gamma$ and $\delta$ are the fold lines.} 
\label{nishi}
\end{figure}

Any quintic equation may be solved by applying the following two-fold operation (see Fig.~\ref{nishi}): given two points $P$ and $Q$ and three lines $\ell$, $m$, $n$, simultaneously fold along a line $\gamma$ to place $P$ onto $m$, and along a line $\delta$ to place $Q$ onto $n$ and to align $\ell$ and $\gamma$ \cite{Nishimura2015}.

\begin{figure}[!htb]
\centering
\includegraphics{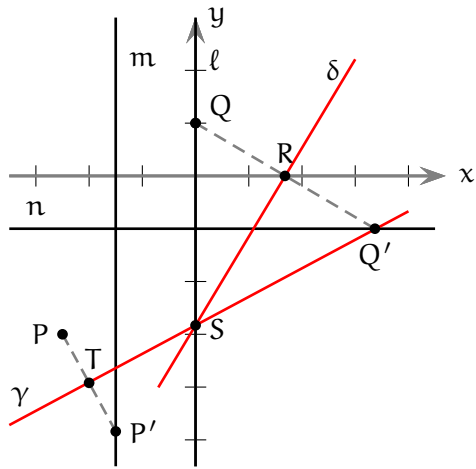}
\caption{Geometric solution of Eq.~(\ref{eqt}), for $t=2\cos(2\pi/11)$.} 
\label{twofold}
\end{figure}

The coordinates of points $P$ and $Q$ and equations of lines $\ell$, $m$, and $n$ are computed from  the coefficients of the quintic equation to solve. In the case of Eq.~(\ref{eqt}), we have  $P(-\frac{5}{2}, -3)$, $Q(0, 1)$, $\ell: x=0$, $m: x=-\frac{3}{2}$, and  $n: y=-1$ (Fig.~\ref{twofold}). Their calculation is omitted here for brevity; complete details of the algorithm may be found in Ref.\ \cite{Nishimura2015}.

Let us demonstrate that the folds solve Eq.~(\ref{eqt}). Folding along line $\delta$ reflects point $Q$ onto $Q'\in n$. Assume that point $Q'$ is located at $(2t, -1)$, where $t$ is a parameter. Then, the slope of segment $\overline{QQ'}$ is $-1/t$, and its midpoint $R$ is at $(t, 0)$.  The fold line $\delta$ is perpendicular to $\overline{QQ'}$ and passes through $R$; therefore, it has an equation
\begin{equation}
y= t(x-t)
\label{delta}
\end{equation}

Folding along line $\gamma$ reflects point $P$ on $P'\in m$. Assume that point $P'$ is located at $(-\frac{3}{2}, 2s)$, where $s$ is another parameter. Then, the slope of segment $\overline{PP'}$ is $2s+3$, and its midpoint $T$ is at $(-2, s-\frac{3}{2})$.  The fold line $\gamma$ is perpendicular to $\overline{PP'}$ and passes through $T$; therefore, it has an equation
\begin{align}
y &= -\frac{x+2}{2s+3}+s-\frac{3}{2}\nonumber\\
 &= -\frac{x}{2s+3}+\frac{2s^2-\frac{13}{2}}{2s+3}
\label{gamma}
\end{align}

Now, the same fold along $\delta$ reflects line $\ell$ over $\gamma$. Let $S$ be the point of intersection of $\delta$ and $\ell$. The $y$-intercept may be obtained by letting $x=0$ in Eq.~(\ref{delta}), which produces $y=-t^2$. Then, the slope of segment $\overline{SQ'}$ is $(t^2-1)/(2t)$. Line $\gamma$ must pass through both $Q'$ and $S$, and therefore it has an equation
\begin{equation}
y= \frac{t^2-1}{2t}x-t^2
\label{gamma2}
\end{equation}
Since Eqs.~(\ref{gamma}) and (\ref{gamma2}) describe the same line, then their respective coefficients must be equal:   
\begin{align}
\frac{-1}{2s+3}&=\frac{t^2-1}{2t}\label{aa}\\
\frac{2s^2-\frac{13}{2}}{2s+3}&=-t^2
\label{chi4}
\end{align}
Solving Eq.~(\ref{aa}) for $s$ and replacing into Eq.~(\ref{chi4}), we obtain Eq.~(\ref{eqt}). Therefore, the $x$-intercept of $\delta$ (i.e., $t$ at point $R$) is a solution of Eq.~(\ref{eqt}). Note that the equation has five possible solutions, and Fig.~\ref{twofold} shows the case of $t=2\cos(2\pi/11)\approx 1.6825\ldots$.  

\section{Folding instructions}

The following diagrams present full instructions for folding the regular hendecagon on a square sheet of paper.
 
Steps (1) to (7) produce points $P$ and $Q$ and lines $\ell$, $m$, $n$ of Fig. \ref{twofold}. The center of the paper is assumed to have coordinates $(0, -1)$, and each side has length 8. In step (1), the vertical and horizontal folds define lines $\ell$ and $n$, respectively. In step (3), the intersection of the fold line with the vertical crease is point $Q$. In step (5), the small crease at the bottom will mark the position of line $m$, after folding the paper backwards in the next step. Finally, step (7) defines point $P$.

 Next, steps (8) and (9) produce the fold lines $\gamma$ and $\delta$, respectively, of Fig.~\ref{twofold}. As a result, point $Q'$ in step (11) is at a distance of $4\cos(2\pi/11)$ from the center of the paper. In the same figure, point $A$ is adopted as one vertex ($z_0$ in Fig.~\ref{hend}) of an hendecagon with radius of 4 units. The fold in step (11) produces a vertical line through $Q'$, and steps (12) and (13) rotate a 4-unit length (from the paper center to point $A$) so as to find the next vertex of the hendecagon (point $D$ in step 14). The remainder steps, (14) to (20)  are used to find the other vertices and folding the sides.

\includegraphics[trim=.6cm 0 0 0]{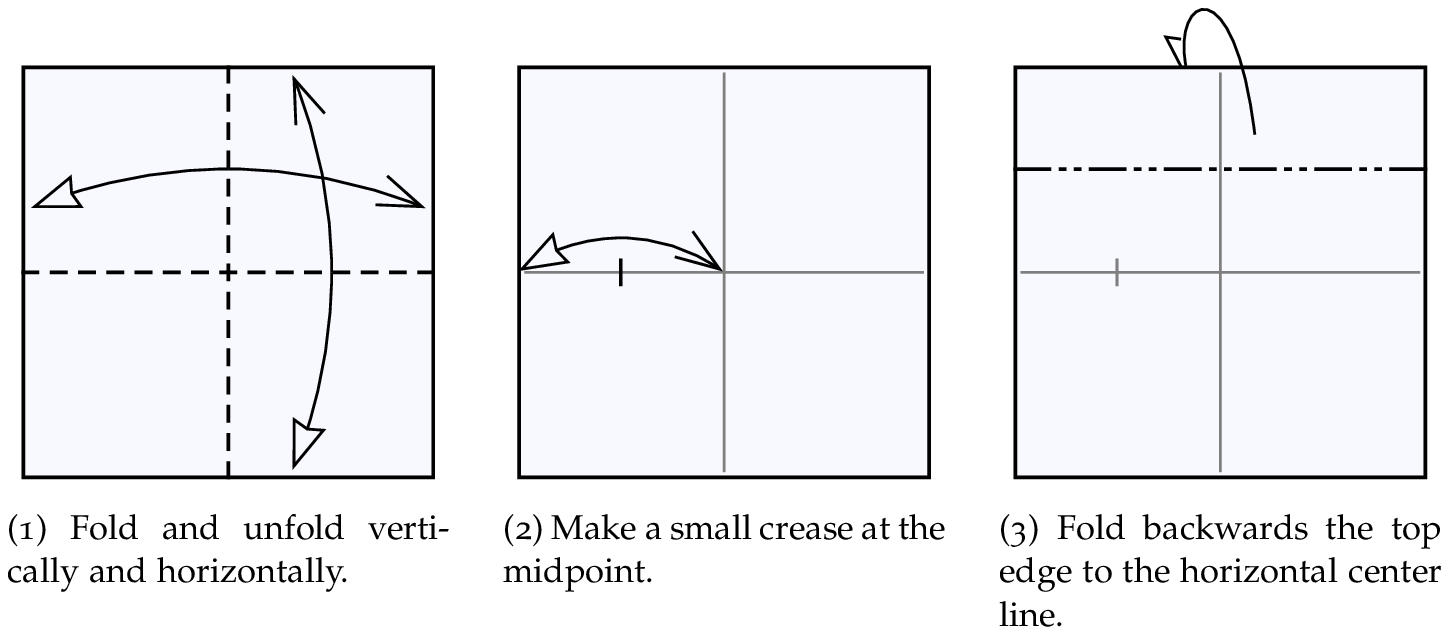}

\bigskip
\bigskip

\includegraphics[trim=.6cm 0 0 0]{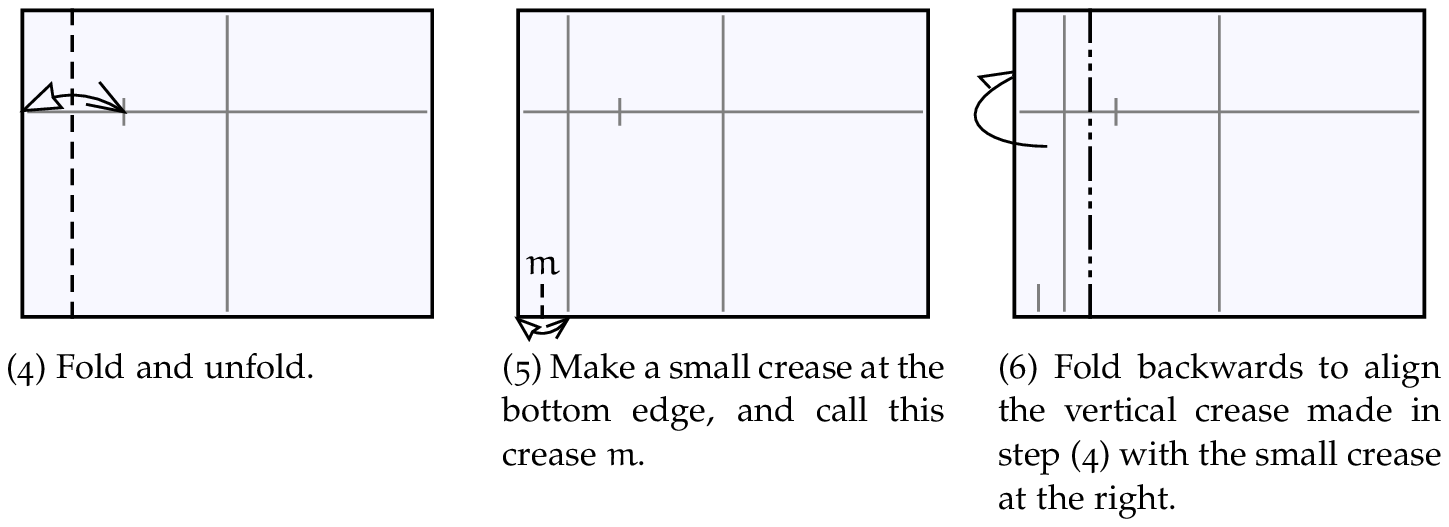}

\bigskip
\bigskip

\includegraphics[trim=.6cm 0 0 0]{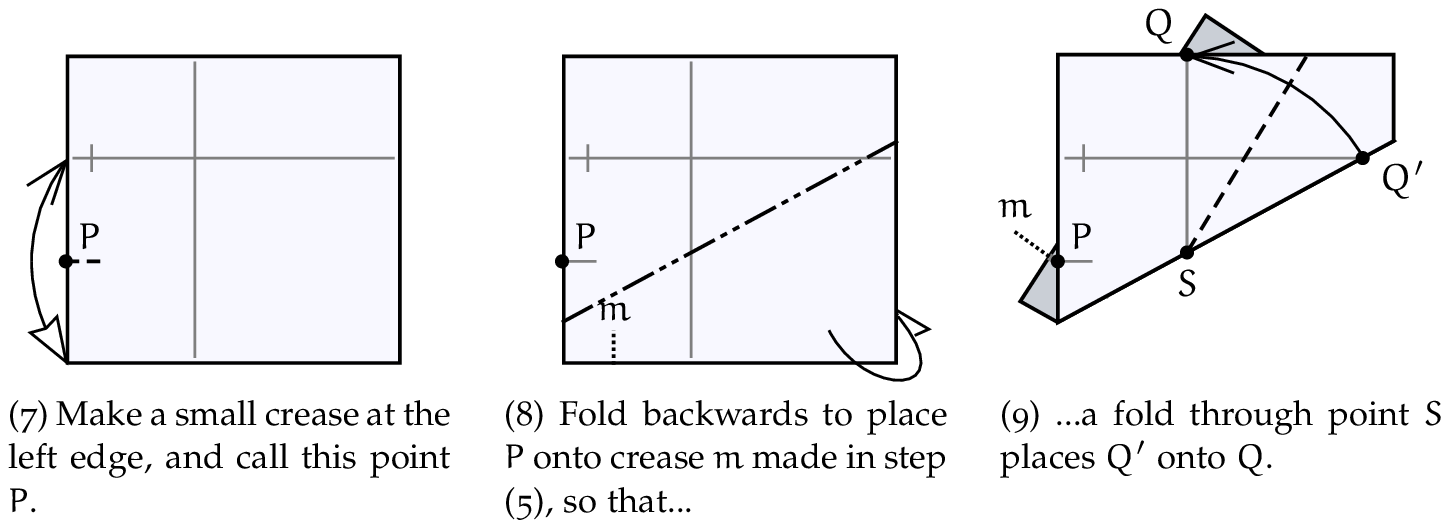}
 
\bigskip
\bigskip

\includegraphics[trim=.6cm 0 0 0]{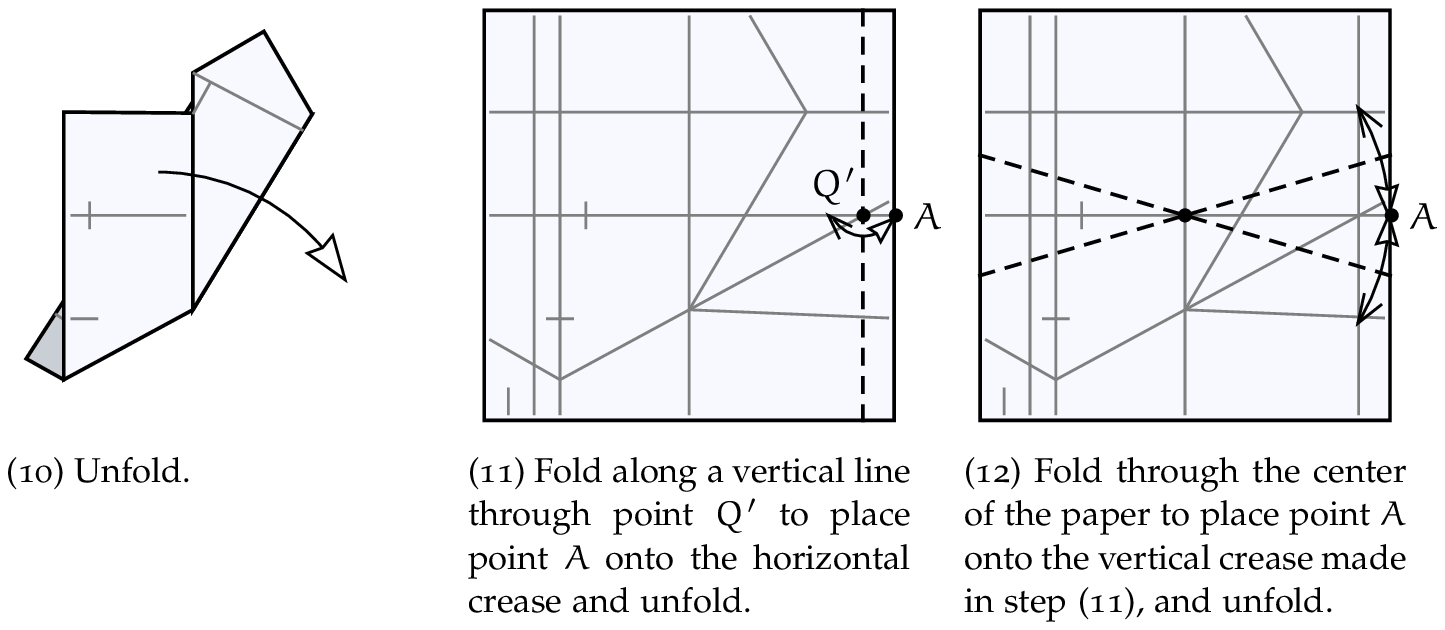}

\bigskip
\bigskip

\includegraphics[trim=.6cm 0 0 0]{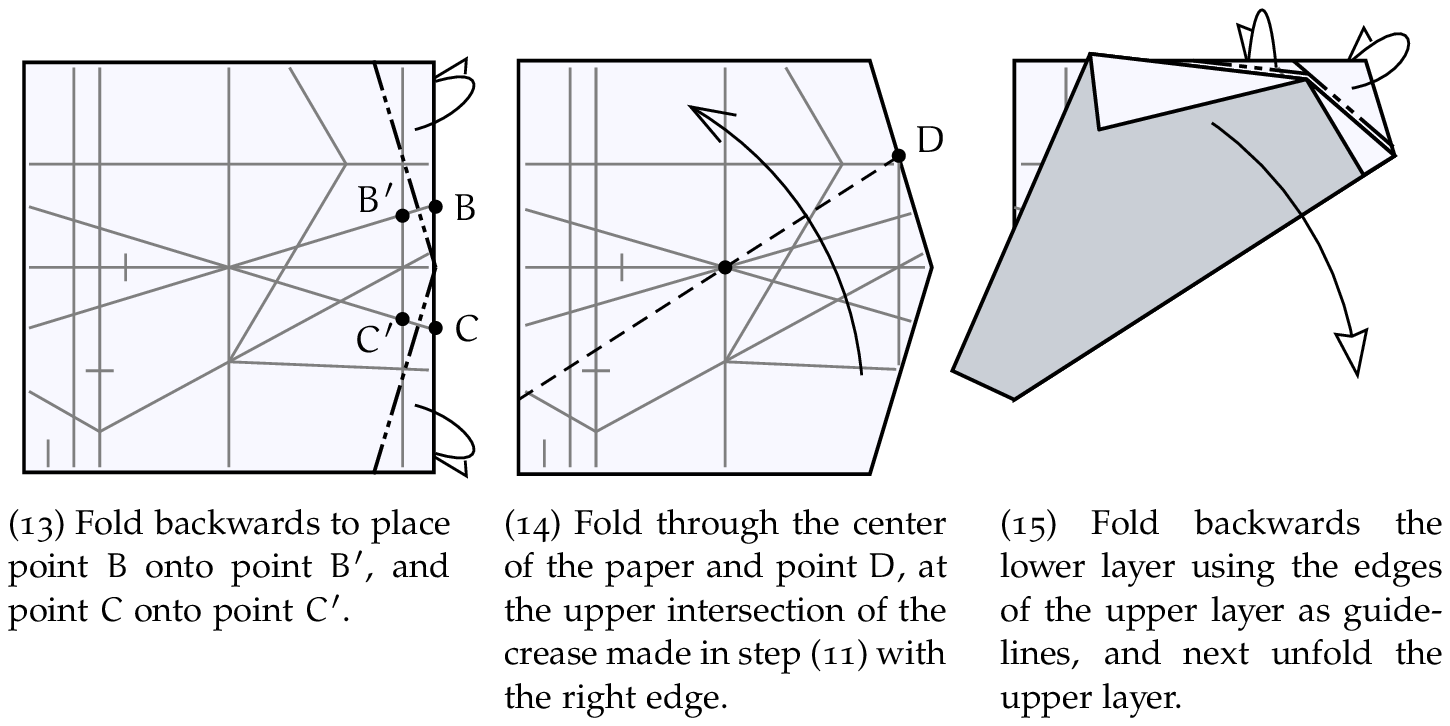}

\bigskip
\bigskip

\includegraphics[trim=.6cm 0 0 0]{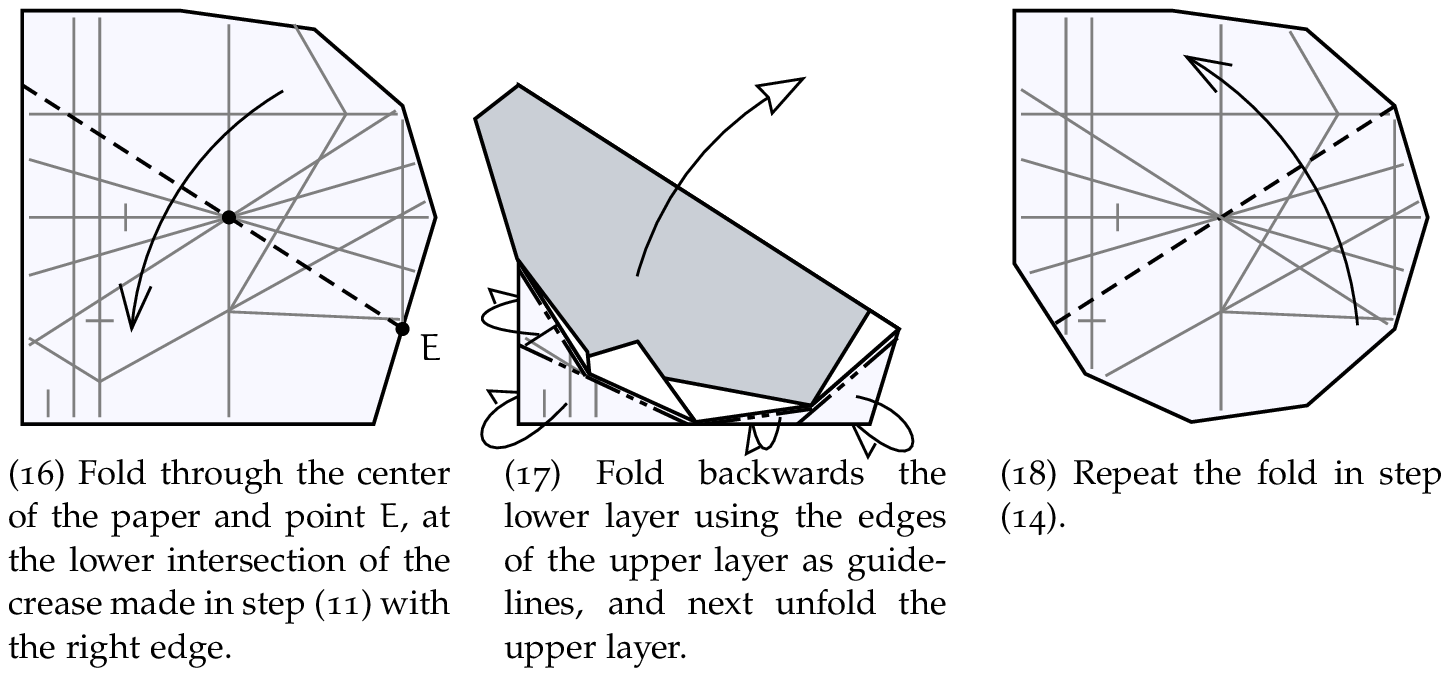}

\bigskip
\bigskip

\includegraphics[trim=.6cm 0 0 0]{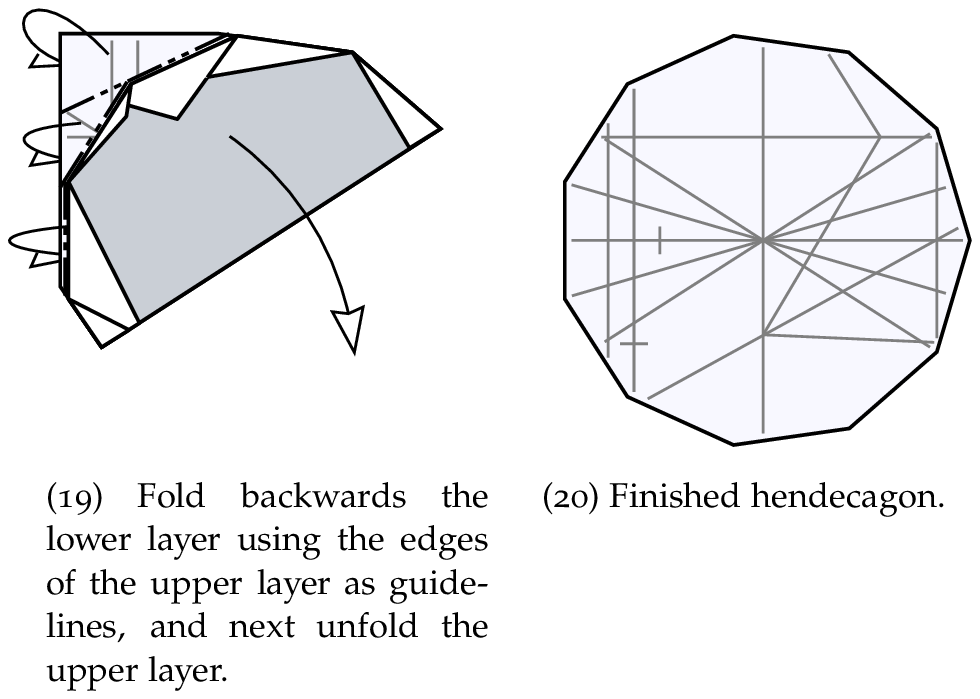}

\end{document}